\documentclass[12pt,twoside,reqno]{amsart}

\newcommand{\titlename}{Itseaffiineista joukoista}

\usepackage[T1]{fontenc}
\usepackage{ae}              % Almost European CM fonts (CM/PS), AMSfonts
\usepackage[finnish]{babel}
\usepackage[latin1]{inputenc}

\usepackage{amsthm}          % theorems
\usepackage{amssymb}         % symbols
\usepackage[dvips]{graphicx} % gfx
\usepackage[bf]{caption2}    % captionstyles
\usepackage{geometry}        % for a good layout on a4
\usepackage{fancyhdr}        % header and footer
\setcaptionwidth{0.425\textwidth}

\long\def\symbolfootnote[#1]#2{\begingroup%
\def\thefootnote{\fnsymbol{footnote}}\footnote[#1]{#2}\endgroup}

\newcommand{\R}{\mathbb{R}}

\newcommand{\N}{\mathbb{N}}

\newcommand{\HH}{\mathcal{H}}
\newcommand{\LL}{\mathcal{L}}
\newcommand{\QQ}{\mathcal{Q}}

\newcommand{\iii}{\mathtt{i}}

\newcommand{\fii}{\varphi}

\DeclareMathOperator{\dimm}{dim_M}

\DeclareMathOperator{\dimh}{dim_H}

\begin{document}

\title{\titlename}

\author{Antti Käenmäki}

%\address{Matematiikan ja tilastotieteen laitos \\
%         PL 35 (MaD) \\
%         40014 Jyväskylän yliopisto}
%\email{antakae@maths.jyu.fi}

%\date{\today}

%\begin{abstract}
%  Arkhimeden artikkeli.
%\end{abstract}

\maketitle

%\section*{\titlename}

\emph{Iteroidulla funktiosysteemillä}\symbolfootnote[0]{\today}
tarkoitetaan äärellistä
kokoelmaa kutistavia kuvauksia $f_1,\ldots,f_k$. Tässä kuvausten
$f_i \colon$ $\R^d \to \R^d$ \emph{kutistavuudella} tarkoitetaan sitä,
että on olemassa $0<\lambda_i<1$ siten, että
\begin{equation}
  \label{eq:kutistus}
  |f_i(x)-f_i(y)| \le \lambda_i|x-y|
\end{equation}
aina kun $x,y\in\R^d$. Hutchinson osoitti vuonna 1981, että jokaiselle
iteroidulle funktiosysteemille on olemassa yksikäsitteinen epätyhjä
kompakti joukko $E \subset \R^d$, jolle
\begin{equation}
  \label{eq:invariantti}
  E = \bigcup_{i=1}^k f_i(E).
\end{equation}
Tuloksen todistus on Banachin kiintopistelauseen elegantti
sovellus. Jos $R>0$ valitaan riittävän suureksi, niin selvästi pätee
$f_i\bigl( B(0,R) \bigr) \subset B(0,R)$, missä $B(x,r)$ on
$x$-keskinen ja $r$-säteinen suljettu pallo. Esimerkiksi valinta $R =
\max_i|f_i(0)|/(1-\max_i\lambda_i)$ riittää. Tällöin on helppo nähdä,
että
\begin{equation}
  \label{eq:peite}
  E = \bigcap_{n=1}^\infty \bigcup_{\iii \in \{ 1,\ldots,k \}^n}
  f_\iii\bigl( B(0,R) \bigr)
\end{equation}
missä $f_\iii = f_{i_1} \circ \cdots \circ f_{i_n}$ aina kun
$\iii = (i_1,\ldots,i_n) \in \{ 1,\ldots,k \}^n$. Tarkempi esitys
aiheesta löytyy esimerkiksi kirjasta \cite[\S 2.2]{Falconer1997}.

Olemme kiinnostuneita joukon $E$ dimensiosta. Käytän tässä yhteydessä
kahta dimensiota, Hausdorffin ja Minkowskin (box counting).
Olkoon $A$ epätyhjä ja rajoitettu $\R^d$:n osajoukko.
Tällöin joukon $A$ \emph{Minkowskin dimensio} on
\begin{align*}
  \dimm(A) &= \inf\{ s : M^s(A) < \infty \} \\
  &= \sup\{ s : M^s(A) > 0 \},
\end{align*}
missä $M^s(A) = \limsup_{\delta \downarrow 0} M_\delta^s(A)$ ja
\begin{gather*}
  M_\delta^s(A) = \inf\{ N\delta^s : A \subset \bigcup_{i=1}^N
  B(x_i,\delta) \\ \text{joillain } x_i \in \R^d \text{ ja } N \in \N \}
\end{gather*}
aina kun $s \ge 0$ ja $\delta > 0$.
Joukon $A$ \emph{Hausdorffin dimensio} taas on
\begin{align*}
  \dimh(A) &= \inf\{ s : \HH^s(A) < \infty \} \\
  &= \sup\{ s : \HH^s(A) > 0 \},
\end{align*}
missä $\HH^s(A) = \lim_{\delta \downarrow 0} \HH^s_\delta(A)$ ja
\begin{gather*}
  \HH^s_\delta(A) = \inf\{ \sum_i r_i^s : A \subset \bigcup_i
  B(x_i,r_i) \\ \text{joillain } x_i \in \R^d \text{ ja } 0 < r_i \le
  \delta \}
\end{gather*}
aina kun $s \ge 0$ ja $\delta > 0$. Molemmat määritelmät antavat
esimerkiksi välille dimension $1$ ja neliölle
dimension $2$. Mutta Hausdorffin ja Minkowskin dimensiot voivat myös
olla ei-kokonaislukuja. Tällaista yleistettyä dimension käsitettä
voidaankin pitää parametrina, joka antaa informaatiota joukon
koosta. Huomaa, että koska selvästi $\HH^s_\delta(A) \le
M^s_\delta(A)$ aina kun $s \ge 0$ ja $\delta > 0$, pätee $\dimh(A) \le
\dimm(A)$. Lisätietoja dimensioista löytyy esimerkiksi kirjoista
\cite[\S 2.1]{Falconer1997}, \cite[\S 2 ja \S 3]{Falconer1990} ja
\cite[\S 4 ja \S 5]{Mattila1995}.

Joukon $E$ dimension määrittäminen täytyy siis aloittaa hyvän peitteen
löytämisellä. Koska molemmat dimensiot ovat skaalausinvariantteja,
voidaan yksinkertaisuuden vuoksi olettaa, että $R = 1$. Tällöin havainnon
\eqref{eq:peite} avulla huomataan, että $f_\iii(E) \subset
B\bigl( f_\iii(0),\lambda_\iii \bigr)$, missä $\lambda_\iii =
\lambda_{i_1} \cdots \lambda_{i_n}$.
Näitä palloja voimme käyttää joukon $E$ peittämiseen. Jotta pääsisimme
käsiksi Minkowskin dimensioon, pitää tarkastella samankokoisia
palloja. Määritellään jokaiselle $0 < \delta < 1$ joukko
\begin{gather*}
  Z(\delta) = \{ \iii : \lambda_\iii < \delta \le
  \lambda_{i_1}\cdots\lambda_{i_{n-1}} \\ \text{jollakin } \iii =
  (i_1,\ldots,i_n) \}.
\end{gather*}
Tällöin jokaiselle jonolle $(i_1,i_2,\ldots)$ löytyy
täsmälleen yksi $n$, jolle $(i_1,\ldots,i_n) \in Z(\delta)$. Edelleen,
jos valitaan $s \ge 0$ siten, että
\begin{equation}
  \label{eq:sim_dim}
  \sum_{i=1}^k \lambda_i^s = 1,
\end{equation}
saadaan induktiolla $\# Z(\delta) (\min_i \lambda_i)^s \delta^s \le
\sum_{\iii \in Z(\delta)} \lambda_\iii^s = 1$ kaikille $0 < \delta <
1$. Näin ollen alkioiden määrä joukossa $Z(\delta)$ on korkeintaan $(\min_i
\lambda_i)^{-s} \delta^{-s}$ ja joukko $E$ voidaan peittää
tällä määrällä $\delta$-säteisiä palloja. Siispä $M_\delta^s(E) \le
(\min_i \lambda_i)^{-s}$ kaikilla $0 < \delta < 1$ ja $\dimm(E) \le s$.

Dimensiolle alarajan löytäminen on paljon vaikeampaa.
Jos ku\-vauk\-set $f_i$ ovat similariteetteja
eli ehdossa \eqref{eq:kutistus} onkin yhtäsuuruus, niin tämä vielä
onnistuu kohtuullisella vaivalla. Hutchinsonin tuloksesta
\eqref{eq:invariantti} seuraa nimittäin suoraan,
että joukko $E$ koostuu pienemmistä ja pienemmistä, $E$:n kanssa
geometrisesti samanlaisista osista. Joukkoa $E$ kutsutaankin tällöin
\emph{itsesimilaariksi joukoksi}.
Jos itsesimilaarin joukon $E$ osat $f_i(E)$ ovat ''pahasti
päällekkäin'', saattaa itsesimilaarin rakenteen
havaitseminen ja hyväksi käyttäminen olla vaikeaa.
Luonnollinen ajatus on tietysti yrittää näyttää, että palloilla
$B\bigl( f_\iii(0),\lambda_\iii \bigr)$ saadaan
muodostettua optimaalisia peitteitä. Tällöin
yhtälön \eqref{eq:sim_dim} määräämä $s$ olisi Hausdorffin dimensiolle
myös alaraja. Hutchinson esittelikin nk.\ avoimen joukon
ehdon, joka takaa palloille riittävän erillisyyden. \emph{Avoimen
  joukon ehdossa} oletetaan, että on olemassa
epätyhjä avoin joukko $V \subset \R^d$, jolle
\begin{equation*}
  V \supset \bigcup_{i=1}^k f_i(V)
\end{equation*}
ja $f_i(V) \cap f_j(V) = \emptyset$ aina kun $i \ne j$.
Hän osoitti, että avoimen joukon ehdon ollessa voimassa on
itsesimilaarin joukon $E$ Hausdorffin mitta positiivinen, $\HH^s(E)
> 0$. Näin ollen $\dimh(E) = \dimm(E) = s$. Todistus tälle on esitetty
kirjassa \cite[Lause 9.3]{Falconer1990}. Schief todisti vuonna
1994, että tulos myös kääntyy: Hausdorffin mitan positiivisuudesta
seuraa avoimen joukon ehto.

Entä jos tarkastellaan yleisempiä kuvauksia?
Oletetaan, että $A_i \in \R^{d \times d}$
on kääntyvä $d \times d$ matriisi, jolla $||A_i|| < 1$, ja
$a_i \in \R^d$ on siirtovektori. Valitsemalla $\lambda_i = ||A_i||$,
affiinit kuvaukset $f_i = A_i + a_i$ muodostavat iteroidun
funktiosysteemin. Joukkoa $E$, jolle \eqref{eq:invariantti} pätee,
sanotaan tällöin \emph{itseaffiiniksi joukoksi}. Yhtälöstä
\eqref{eq:sim_dim} saadaan tietysti yläraja itseaffiinin
joukon dimensiolle. Mutta koska
$f_\iii\bigl( B(0,1) \bigr) = A_\iii\bigl( B(0,1) \bigr) + a_\iii$
jollakin $a_\iii \in \R^d$, voidaan joukko $E$ peittää pallojen
sijasta ellipseillä ja siten mahdollisesti parantaa ylärajaa.
\begin{figure}
  \centering
  \includegraphics[width=0.35\textwidth]{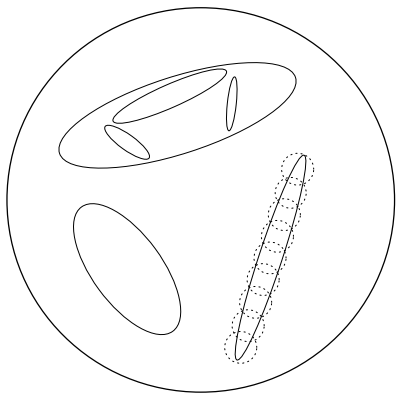}
  \caption{Itseaffiinia joukkoa voidaan peittää ellipseillä, joita
    taas voidaan peittää pienemmillä palloilla.}
  \label{fig:ellipses}
\end{figure}
Tässä $A_\iii = A_{i_1} \cdots A_{i_n}$ aina kun $\iii = (i_1,\ldots,i_n)$.
Olkoon ellipsin $f_\iii\bigl( B(0,1) \bigr)$ puoliakselien pituudet $1
> \alpha_1(A_\iii) \ge \cdots \ge \alpha_d(A_\iii) > 0$.
Nämä pituudet saadaan laskettua ottamalla positiivisesti definiitin
matriisin $A_\iii^T A_\iii$ ominaisarvoista neliöjuuret. Koska
esimerkiksi tasossa $\alpha_2(A_\iii)$-säteisiä palloja tarvitaan
ellipsin $f_\iii\bigl( B(0,1) \bigr)$ peittämiseen oleellisesti noin
$\alpha_1(A_\iii)/\alpha_2(A_\iii)$ kappaletta, on
\begin{equation*}
  \HH^s(E) \lesssim \lim_{n \to \infty} \sum_{\iii \in \{ 1,\ldots,k
    \}^n} \alpha_1(A_\iii)\alpha_2(A_\iii)^{s-1}
\end{equation*}
ja siten ylärajan löytämiseksi Hausdorffin dimensiolle riittää
tarkastella kuinka yo.\ summa käyttäytyy eri $s$:n arvoilla. Tämän
esimerkin motivoimana määrittelemme jokaiselle $0 \le
s < d$
\begin{equation*}
  \fii^s(A_\iii) = \alpha_1(A_\iii) \cdots \alpha_{l+1}(A_\iii)^{s-l},
\end{equation*}
missä $l$ on $s$:n kokonaislukuosa, ja
\begin{equation*}
  P(s) = \lim_{n \to \infty} \tfrac{1}{n} \log \sum_{\iii \in \{
    1,\ldots,k \}^n} \fii^s(A_\iii).
\end{equation*}
Raja-arvo edellä on olemassa, koska kyseessä oleva jono on
subadditiivinen (tämän todistamiseksi tarvitaan
multilineaarialgebraa). Nyt nähdään helposti, 
että $\sum \fii^s(A_\iii) \asymp e^{nP(s)}$ ja että $P$:llä on
yksikäsitteinen nollakohta, jota kutsutaan
\emph{singulaaridimensioksi}.
Vuonna 1988 Fal\-coner näytti, että itseaffiinilla joukolla
singulaaridimensio on aina yläraja Minkowskin dimensiolle. Huomaa
myös, että itsesimilaarissa tilanteessa yhtälön \eqref{eq:sim_dim}
määräämä $s$ on sama kuin singulaaridimensio.

Koska itseaffiinilla joukolla singulaaridimensio on alaraja
Hausdorffin tai Minkowskin dimensiolle? Pitäisi siis taas pystyä
muodostamaan optimaalisia peitteitä. Nyt pahan päällekkäisyyden
lisäk\-si meillä on
kaksi uutta tilannetta, mitkä saattavat aiheuttaa ongelmia.
Ellipsit voivat olla sijoittuneet niin, että yhden vähän isomman
pallon käyttäminen peitteenä antaisi paremman lopputuloksen kuin
kyseessä olevien ellipsien peittäminen pienemmillä palloilla. Lisäksi
itseaffiinia joukkoa pitäisi olla ''riittävästi'' ellipsin pisimmän
akselin suuntaisesti, jotta ellipsiä peitettäessä pienempiä palloja ei
käytettäisi turhaan.
\begin{figure}[b]
  \centering
  \includegraphics[width=0.45\textwidth]{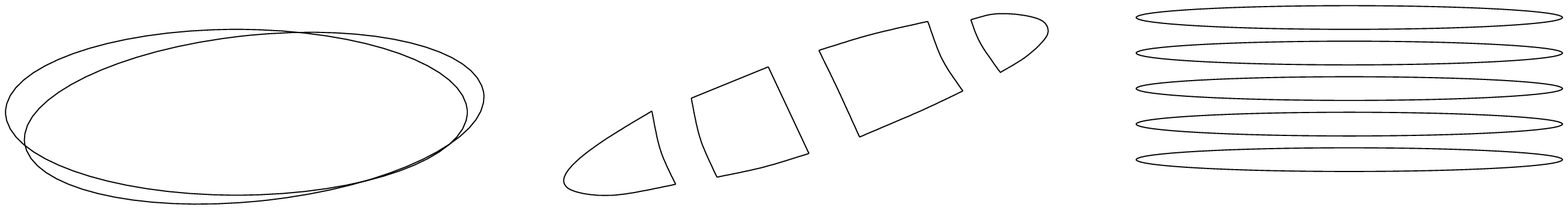}
  \caption{Ellipsien paha päällekkäisyys, ''reikäisyys'' ja huono
    sijoittuminen saattavat pilata peitteen optimaalisuuden.}
  \label{fig:drop}
\end{figure}
Kuva \ref{fig:drop} havainnollistaa näitä tilanteita. Huomaa myös,
että peittävät ellipsit saattavat olla päällekkäin sillä tavalla,
että pidempien akseleiden välinen kulma on suuri. Tämä ei välttämättä
ole optimaalisuuden kannalta huono tilanne.

Bedford ja McMullen (toisistaan riippumatta) laskivat vuonna 1984 ns.\
itseaffiineille matoille Hausdorffin ja Minkowskin
dimensiot. \emph{Itseaffiinilla matolla}
tarkoitetaan joukkoa, joka saadaan seuraavanlaisella konstruktiolla:
Jaetaan neliö $p$ moneen samankokoiseen sarakkeeseen ja
$q$ moneen samankokoiseen riviin ja näin saaduista suorakaiteista
valitaan osa. Jokaisessa valitussa suorakaiteessa tehdään samanlainen
jako ja valitaan näissä vastaavat suorakaiteet kuin alussakin. Näin
jatkamalla jäljelle jää joukko, jota sanotaan itseaffiiniksi
matoksi. Itseaffiini matto on selvästi itseaffiini joukko, sillä
affiineiksi kuvauksiksi voidaan valita kuvaukset, jotka vievät neliön
valituiksi suorakaiteiksi.
\begin{figure}
  \centering
  \includegraphics[width=0.475\textwidth]{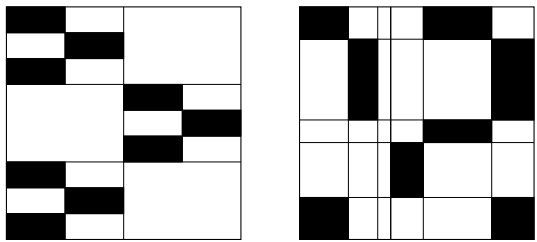}
  \caption{Bedford-McMullenin ja Bara\'nskin itseaffiinit matot.}
  \label{fig:carpets}
\end{figure}
Kuva \ref{fig:carpets} havainnollistaa tilannetta. Osoittautuu, että
sopivilla suorakaiteiden valinnoilla tällaiselle
matolle $E$ saadaan $\dimh(E) < \dimm(E)$ $< s$, missä $s$ on
singulaaridimensio. Dimensio siis riippuu siirtovektoreiden
valinnasta! Katso \cite[Esimerkki 9.11]{Falconer1990}. Bedfordin ja
McMullenin konstruktiota on sittemmin
yleistetty monella eri tapaa. Esimerkiksi vuonna 2007 Bara\'nski
määritti hyvin yleisille itseaffiineille matoille Hausdorffin ja
Minkowskin dimensiot. Hänen konstruktiossaan neliö voidaan jakaa
suorakaiteisiin hyvin vapaasti: mikä tahansa ositus, joka saadaan
äärellisillä määrillä vaakasuoria ja pystysuoria viivoja kelpaa.

Näiden esimerkkien valossa näyttäisi kovasti siltä, että useimmiten
itseaffiinin joukon dimensio eroaisi
singulaaridimensiosta. Kuitenkin Falconer o\-soit\-ti vuonna 1988,
että jos $s$ on singulaaridimensio ja
$||A_i|| < \tfrac13$ kaikilla $i$, niin $\dimh(E_a) = s$
$\LL^{dk}$-melkein jokaisella siirtovektoreiden $a = (a_i,\ldots,a_k)$
valinnalla. Tässä $\LL^{dk}$ on $dk$-ulotteinen Lebesguen mitta ja
merkintä $E_a$ tarkoittaa sitä, että itseaffiini joukko $E$
riippuu siirtovektoreista, mutta matriisit $A_i$ ovat
kiinnitetty. Tulos on hämmästyttävä: se osoittaa, että tyypillisellä
itseaffiinilla joukolla ellipseillä muodostetut peitteet ovat
optimaalisia. Edgar näytti samana vuonna esimerkillä, että tulos ei
pidä paikkaansa, jos yhdelläkin matriisilla on $||A_i||$ $>
\tfrac12$. Solomyak taas osoitti vuonna 1998, että tuloksessa matriisien 
normin yläraja $\tfrac13$ voidaan korvata $\tfrac12$:lla.
Vuonna 2007 Jordan, Pollicott ja Simon todistivat satunnaisen version
Falconerin tuloksesta. Olettamatta matriisien normeille ylärajaa he
osoittivat, että todennäköisyydellä $1$ satunnaisen itseaffiinin
joukon Hausdorffin dimensio on sama kuin singulaaridimensio.
Tässä satunnainen itseaffiini joukko muodostetaan kuten kohdassa
\eqref{eq:peite} paitsi, että konstruktion joka vaiheessa affiinin
kuvauksen siirtovektoriin lisätään satunnainen virhe.
Mainittakoon vielä, että Falconer ja Miao ovat (vielä
julkaisemattomassa artikkelissaan) tutkineet niiden siirtovektoreiden
$a$, joilla $\dimh(E_a) < t$ annetulla $t$, muodostaman joukon
Hausdorffin dimensiota. He osoittivat, että tämä dimensio on korkeintaan
$dk - c(s-t)$, missä $s$ on singulaaridimensio ja $c$ jokin
positiivinen vakio.

\begin{figure}
  \centering
  \includegraphics[width=0.475\textwidth]{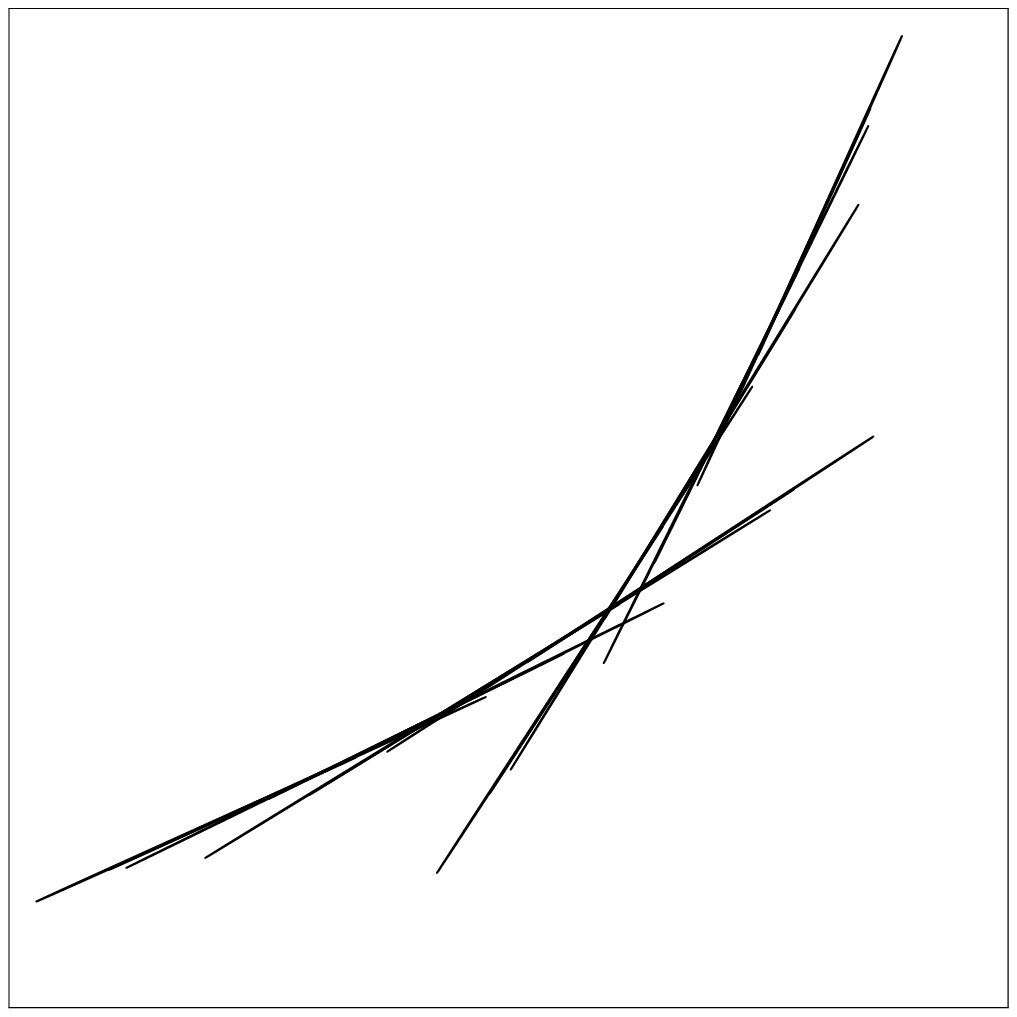}
  \caption{Itseaffiini joukko $E$, jolla $\dimm(E) \ge 1$.}
  \label{fig:kakeya}
\end{figure}

Koska Falconerin tulos ei kerro onko annetun itseaffiinin joukon
dimensio sama kuin singulaaridimensio, on mielenkiintoista yrittää
löytää ehtoja, joiden voimassa ollessa näin kävisi. Falconer o\-soit\-ti
vuonna 1992, että jos $s$ on singulaaridimensio, avoimen joukon ehto
on voimassa yhtenäiselle joukolle $V \subset \R^d$ ja joukon $E$
projektiolla mille tahansa $(d-1)$-ulotteiselle aliavaruudelle on
positiivinen Lebesguen mitta, niin $\dimm(E) = s$. Huomaa, että
täytyy olla $\dimh(E) \ge d-1$, jotta projektioehto voisi olla
voimassa. Ehdolla varmistetaankin, että itseaffiinia joukkoa on
''riittävästi'' peittävien ellipsien pisimmän akselin suuntaisesti.
Vuonna 1995 Hueter ja Lalley esittelivät tason itseaffiineille
joukoille ehdot, joiden voimassa ollessa $\dimh(E) = s < 1$, missä $s$
on singulaaridimensio. He olettivat, että peittävien
ellipsien muoto on rajoitettu, $\alpha_1(A_i)^2 < \alpha_2(A_i)$, ja
että ne ovat riittävästi erillään. Lisäksi he vaativat
lineaarikuvausten $A_i$ kuvaavan tason ensimmäisen neljänneksen
$\QQ_2$ osakseen sillä tavalla, että $A_i(\QQ_2) \cap A_j(\QQ_2) =
\emptyset$ aina kun $i \ne j$. Koska $s<1$, riippuu dimensio Hueterin
ja Lalleyn tapauksessa vain ellipsien pisimmän puoliakselin
pituudesta. Shmerkin ja allekirjoittanut ovat (vielä
julkaisemattomassa artikkelissa) esitelleet tasossa ehdot, joiden
voimassa ollessa $\dimm(E) = s \ge 1$, missä $s$ on
singulaaridimensio. Ehdoissa oletetaan samantapainen vaatimus
koskien lineaarikuvauksia $A_i$ ja neljännestä $\QQ_2$ kuin Hueterilla
ja Lalleylla sekä joukon $E$ projektio mille tahansa suoralle, jolla on
positiivinen kulmakerroin, oletetaan positiivimittaiseksi. Koska
nyt $s \ge 1$, dimensio riippuu ellipsien molempien puoliakselien
pituuksista. Lisäksi huomion arvoista on, että matriisien normeille ei
ehdoissa vaadita ylärajaa ja että peittävillä ellipseillä saa olla
päällekkäisyyttä. Tuloksessa on käytetty hyväksi Kakeya-jouk\-ko\-jen
teoriaa. Lisätietoja Kakeya-jou\-kois\-ta löytyy esimerkiksi Arkhimeden
artikkelista \cite{Mattila2004}.
Kuvassa \ref{fig:kakeya} on esimerkki nämä ehdot
täyttävästä itseaffiinista joukosta.

Itseaffiinin joukon dimension määrittäminen on siis hankalaa. Jopa
yksinkertaisen oloisissa tilanteissa (esim.\ itseaffiinit matot) dimension
selville saaminen on hyvin vaikeaa. Jos kuitenkin saadaan osoitettua
yleinen tulos (esim.\ Falconerin tulos), on tulos geneerinen eli se ei sano
mitään annetulle itseaffiinille joukolle. Riittävien ehtojen
löytämiseksi täytyy usein tehdä hyvin rajaavia oletuksia (esim.\
Hueterin ja Lalleyn tulos). Hankalaa on myös singulaaridimension
laskeminen. Falconer ja Miao löysivät vuonna 2007
singulaaridimensiolle suljetun kaavan olettamalla matriisit $A_i$
yläkolmiomatriiseiksi.

Mainitaan vielä lopuksi, että itseaffiineilla joukoilla on käytännön
sovelluksia kuvankäsittelyssä. Kirjoilla \cite[\S 9.5]{Falconer1990}
ja \cite[\S 9.8]{Barnsley1988} pääsee aiheessa alkuun.

\bibliographystyle{abbrv}
\bibliography{itseaffiineista.bib}

\end{document}